\documentclass[12pt,reqno]{amsart}
\usepackage[all]{xy}
\usepackage{mathrsfs}
\usepackage{amsmath, amsthm, amssymb, amscd, amsgen,amsxtra,amsfonts, amsbsy}
\allowdisplaybreaks
\usepackage[all]{xy}
\usepackage{verbatim}

\usepackage{tikz}
\usepackage{pgfplots}
\usepackage{tikz-cd}
\usepackage{etex}
\usepackage{calligra,mathrsfs}
\usepackage{color}
\definecolor{shadecolor}{gray}{0.875}
\definecolor{dblue}{rgb}{1,0,0}
\usepackage[colorlinks=true, linkcolor=dblue, citecolor=dblue, filecolor = dblue, menucolor = dblue, urlcolor = dblue]{hyperref}
\usepackage{mathtools}
\usepackage{extarrows}
\usepackage{ifthen}
\usepackage{bm}

\usepackage{graphicx}
\usepackage{caption,rotating}
\usepackage{color}
\usepackage{subcaption}

\usepackage{enumitem}

\numberwithin{equation}{section}

\usepackage[toc,page]{appendix}

\renewcommand{\deg}{\operatorname{deg}}

\newcommand{\N}{\mathbb{N}}

\newcommand{\Q}{\mathbb{Q}}
\newcommand{\C}{\mathbb{C}}
\newtheorem{theorem}{Theorem}[section]
\newtheorem{lemma}[theorem]{Lemma}
\newtheorem{corollary}[theorem]{Corollary}
\newtheorem{conjecture}[theorem]{Conjecture}
\newtheorem{proposition}[theorem]{Proposition}
\newtheorem*{theorem*}{Theorem}
\newtheorem*{corollary*}{Corollary}
\theoremstyle{definition}
\newtheorem{defn}[theorem]{Definition}

\newtheorem{remark}[theorem]{Remark}

\numberwithin{equation}{section}

\author{Houari Benammar Ammar}
\AtEndDocument{\bigskip{\footnotesize%
  \textsc{Department of Mathematical and Statistical Sciences, University of Alberta, Edmonton, Alberta T6G 2G1, Canada} \par  
  \textit{E-mail address :} \texttt{hbenamma@ualberta.ca} 
}}
\large
\title{Note on Iitaka Conjecture $C_{n,m}$}
\begin{document}

\normalsize
\maketitle
\begin{abstract}
Let \( f: X \to Y \) be an algebraic fiber space, where \( X \) and \( Y \) are smooth projective varieties of dimensions \( n \) and \( m \), respectively.  
In \cite{Caopaun}, Cao and Păun proved \( C_{n,m} \) when \( Y \) has maximal Albanese dimension.  
In this paper, we prove \( C_{n,m} \) in the case where the Albanese dimension $\alpha(Y)$ of \( Y \) satisfies \( \alpha(Y) \geq m-2 \).
\end{abstract}
\section{Introduction}
In this paper, we work over $\C$. An algebraic fiber space is a surjective morphism $f: X \to Y$ of
smooth varieties, with connected fibers. A fundamental problem
in birational geometry is the following Iitaka conjecture, denoted by $C_{n,m}$. 
\begin{conjecture}\label{Conj1}
Let $f : X \to Y$ be an algebraic fiber space, where $X$ and $Y$ are smooth projective varieties of dimension $n$ and $m$,
respectively, and let $F$ be a general fiber of $f$. Then
\[
\kappa(X) \geq \kappa(F) + \kappa(Y).
\]
\end{conjecture}
Kawamata~\cite{Kawamatacurves} proved \( C_{n,1} \) and also established the conjecture in the case where the general fiber \( F \) admits a good minimal model~\cite{kawamataminimalmodel}. In particular, his result covers Kollár's result~\cite{kollargeneraltype} when \( F \) is of general type, as well as Fujino's results~\cite{fujinomaximal} when \( F \) has maximal Albanese dimension. 
Indeed, by~\cite{bchm} and~\cite{fujinogood}, varieties of general type and varieties with maximal Albanese dimension admit good minimal models, respectively. Viehweg \cite{viehweglitaka} proved $C_{n,m}$ when $Y$ is of general type. Moreover, Birkar~\cite{Birkar_2009} proved \( C_{n,m} \) for \( n \leq 6 \), and Cao~\cite{Caolitaka2} proved \( C_{n,2} \).  It is also worth mentioning that Campana \cite{campanaorbifold} and Lu \cite{Steven} proved several results for the orbifold $C_{n,m}$ conjecture.
Here we are particularly interested in the case where \( Y \) is irregular (and hence \( X \) is irregular as well). 
Recall that Cao and Păun~\cite{Caopaun} settled \( C_{n,m} \) when \( Y \) has maximal Albanese dimension (see \cite{HSP} for a nice exposition). 
In this paper, we prove the following theorem.
\begin{theorem}\label{maintheorem}
Iitaka conjecture $C_{n,m}$ holds when $\alpha(Y) \geq m-2$.
\end{theorem}
\subsection*{Acknowledgments}
I am grateful to Kenta Hashizume, who found a minor mistake in the previous version, and I would like to think him for informing me. I would also like to thank Xi Chen for his remarks, insightful questions, and encouragement. I decided to write this paper when Mihai Păun visited Steven Lu in Montreal in August 2025. I benefited from discussions with both of them, as well as from conversations with Emanuele Ronda on the subject. A special thanks to Steven for his engagement and inestimable help.

\section{Preliminaries}
\begin{defn}
     We define the Albanese dimension $\alpha(X)$ of an irregular variety $X$ by 
    \[
    \alpha(X):= \dim \operatorname{alb}_{X}(X).
    \]
    Here, $\operatorname{alb}_{X}(X)$ is the image of the Albanese map $\operatorname{alb}_{X}: X \to \operatorname{alb}_{X}(X) \subseteq \operatorname{Alb}(X)$, where $\operatorname{Alb}(X)$ is the Albanese variety.
    \par We call $\alpha_{f}(X): = \dim X - \alpha(X)$ the Albanese fiber dimension of $X$. By definition $ 0\leq \alpha(X) \leq \dim X$ and $0 \leq \alpha_{f}(X) \leq \dim X$.
\end{defn}
Recall the following proposition due to Fujino, which will be indispensable for the proof of our main Theorem~\ref{maintheorem}.
\begin{proposition}[{\cite[Proposition 2.5]{fujinoremarks}}]\label{propo}
Let $V$ be a non-singular projective variety and $W$ a closed subvariety of $V$. If $W$ is general in $V$, then 
\[
\alpha_{f}(W) \leq \dim W - \dim \operatorname{alb}_{V}(W) \leq \alpha_{f}(V).
\]
\end{proposition}
The next Lemma is due to Kawamata as mentioned by Fujino \cite{fujinomaximal}.
\begin{lemma}[{\cite[Lemma 4.3]{fujinomaximal}}]\label{lemma}
Let $f : X \to Y$ and $g : X \to Z$ be surjective morphisms between
normal projective varieties with connected fibers. Then there exist a general subvariety $\Gamma$ of $Z$ , a projective variety $\widehat{Y}$ and morphisms $\gamma: \widehat{Y} \to \Gamma  , \beta : \widehat{Y} \to Y$ such that
$\gamma$ is a surjective morphism with connected fibers, and $\gamma^{-1}(z) \simeq f(g^{-1}(z))$ for $z \in \Gamma$, and $\beta$ is generically finite.
\end{lemma}
\par During the proof of Theorem \ref{maintheorem}, especially in Step 4, we need to apply the Fujino-Mori canonical bundle formula.
\begin{theorem}[{\cite[Theorem 4.5]{fujinomori}}]\label{fujinomoricbf}
Let $f: X \to Y$ be an algebraic fiber space, and let $F$ be a general fiber of $f$ with $\kappa(F) =0$. Then, there exists a commutative diagram:
\[
\begin{tikzcd}
\Tilde{X} \arrow[r, "\Tilde{f}"] \arrow[d, "\sigma"] & \Tilde{Y} \arrow[d, "\tau"] \\
X \arrow[r, "f"] & Y
\end{tikzcd}
\]    
with the following properties:
\begin{itemize}
    \item[1)]$\Tilde{f}$ is an algebraic fiber space of smooth varieties, $\sigma$ and $\tau$ are birational.
\item[2)] There is a $\Q$-effective divisor $B$ on $\Tilde{Y}$, a nef $\Q$-divisor $L$ on $\Tilde{Y}$ and a $\Q$-divisor $R = R^{+}-R^{-}$ on $\Tilde{X}$ decomposed into its positive and negative parts such that
\[
pK_{\Tilde{X}} = p\Tilde{f}^{*}(K_{\Tilde{Y}} + B + L) + R,
\]
for a sufficiently integer $p \in \N$.
\item[3)] $\Tilde{f}_{*}\mathcal{O}_{\Tilde{X}}(iR^{+}) = \mathcal{O}_{\Tilde{Y}}$ for any $i \in \N$.
\item[4)] The divisor $R^{-}$ is exceptional $/X$ and the codimension of $\Tilde{f}(\text{Supp} R^{-})$ in $\Tilde{Y}$ is $\geq 2$.
\end{itemize}

\end{theorem}
\par Independently of our proof, we give an application of \cite[Theorem 12.1]{schnellcampanaconj} below, and provide a quick proof of Theorem \ref{maintheorem} in the case where the Iitaka variety of $Y$ is not uniruled.
\begin{theorem}[{\cite[Theorem 12.1]{schnellcampanaconj}}]\label{schnellthm}
Let $f: X \to Y$ be an algebraic fiber space with $\kappa(F) \geq 0$. Let $H$
be an ample divisor on $Y$. If $mK_{X} - f^{*}H$ is pseudo-eﬀective for some $m\geq 1$, then
$m K_{X} - f^{*}H$ becomes eﬀective for $m$ suﬃciently large and divisible, provided that $Y$ is not uniruled.
\end{theorem}


\section{Iitaka conjecture with big Albanese dimension}
\begin{proof}[Proof of Theorem \ref{maintheorem}] Let $f: X \to Y$ be the given algebraic fiber space. We divide the proofs is several steps.
\par Step 1. If $\kappa(F) = -\infty$ or $\kappa(Y) = -\infty$, then the inequality is satisfied. From now, we assume that $\kappa(F) \geq 0$ and $\kappa(Y) \geq 0$.
\par Step 2. We use the Albanese map $\operatorname{alb}_{Y}: Y \to \operatorname{alb}_{Y}(Y) \subseteq \operatorname{Alb}(Y)$, and obtain the following diagram
\[
\begin{tikzcd}
X \arrow[r, "f"] \arrow[rd, "h"] & Y \arrow[d, "\operatorname{alb}_{Y}"] \\
 & \operatorname{alb}_{Y}(Y)
\end{tikzcd}
\]
Up to birational modification and Stein factorization, we assume that $\operatorname{alb}_{Y}(Y)$ is smooth, $\operatorname{alb}_{Y}$ is an algebraic fiber space and $G$ its general fiber. We denote by $h$ the composition morphism $h:=\operatorname{alb}_{Y} \circ f$, and by $H$ its general fiber. Here $\operatorname{alb}_{Y}(Y)$ has maximal Albanese dimension. Thus, we apply \cite{Caopaun} to the morphism $h$, and obtain
\begin{equation}\label{eq1}
    \kappa(X) \geq \kappa(H) + \kappa(\operatorname{alb}_{Y}(Y)).
\end{equation}
We define the restriction $f_{|_{H}} : H \to G$. Here $F$ is a general fiber of $f_{|_{H}}$. By assumption $\alpha(Y) = m-1$ or $\alpha(Y) = m-2$, which implies that $\dim G = \alpha_{f}(Y) =1$ or $\dim G= \alpha_{f}(Y) =2$, respectively. Thus, we apply \cite{Caolitaka2} and \cite{Kawamatacurves} to obtain
\begin{equation}\label{eq2}
    \kappa (H) \geq \kappa(F) + \kappa(G).
\end{equation}
Combining inequalities (\ref{eq1}) and (\ref{eq2}), we obtain 
\begin{equation}\label{eq3}
    \kappa(X) \geq \kappa(F) + \kappa(G) + \kappa(\operatorname{alb_{Y}(Y)}).
\end{equation}
\par Step 3. $\kappa(Y)=0$: Recall that \(\operatorname{alb}_{Y}(Y)\) has maximal Albanese dimension. Hence, \(\kappa(\operatorname{alb}_{Y}(Y)) \geq 0\) by \cite{Kawamatacarecterization}. Applying \cite{Caopaun} to \(\operatorname{alb}_{Y}\), we obtain
\begin{equation}\label{eq4}
0 = \kappa(Y) \geq \kappa(G) + \kappa(\operatorname{alb}_{Y}(Y)).
\end{equation}
Since \(\kappa(Y) = 0\), it follows that \(\kappa(G) \geq 0\). By inequality~\eqref{eq4}, we deduce that \(\kappa(G) = \kappa(\operatorname{alb}_{Y}(Y)) = 0\). Therefore, by inequality~\eqref{eq3}, we obtain
\[
\kappa(X) \geq \kappa(F),
\]
and the theorem is proved.
\par Step 4. $\kappa(X) =0 \implies \kappa(Y) =0$: By assumption, $\kappa(X) =0$. Thus, by inequality (\ref{eq3}), we have $\kappa(F)= \kappa(G) = \kappa (\operatorname{alb}_{Y}(Y)) =0$. Since the variety $\operatorname{alb}_{Y}(Y)$ is of maximal Albanese dimension, it is birational to an abelian variety, and hence we may assume that it is an abelian variety $A$ (i.e. $\operatorname{alb}_{Y}(Y) =A$). By \cite{birkarchen}, and \cite{Huzhengyu}, we can run the minimal model program (MMP) for $Y$ and we obtain a good MMP $Y^{'}$ of $Y$. Indeed, the variety $Y$ is fibered over an abelian variety, with abundant general fiber $G$, which is either of dimension $2$ with $\kappa(G)= 0$, or is an elliptic curve. Without loss of generality, we write $Y^{'}= Y$, so that $K_{Y}$ is semi-ample. Applying Theorem \ref{fujinomoricbf} to $f$, there exists a commutative diagram
\[
\begin{tikzcd}
\Tilde{X} \arrow[r, "\Tilde{f}"] \arrow[d, "\sigma"] & \Tilde{Y} \arrow[d, "\tau"] \\
X \arrow[r, "f"] \arrow[rd, "h"] & Y \arrow[d, "\operatorname{alb}_{Y}"] \\
 & A
\end{tikzcd}
\] 
with the following properties:
\begin{itemize}
    \item[1)]$\Tilde{f}$ is an algebraic fiber space of smooth varieties, $\sigma$ and $\tau$ are birational.
\item[2)] There is a $\Q$-effective divisor $B$ on $\Tilde{Y}$, a nef $\Q$-divisor $L$ on $\Tilde{Y}$ and a $\Q$-divisor $R = R^{+}-R^{-}$ on $\Tilde{X}$ decomposed into its positive and negative parts such that
\begin{equation}\label{eq5}
pK_{\Tilde{X}} = p\Tilde{f}^{*}(K_{\Tilde{Y}} + B + L) + R.
\end{equation}
\item[3)] $\Tilde{f}_{*}\mathcal{O}_{\Tilde{X}}(iR^{+}) = \mathcal{O}_{\Tilde{Y}}$ for any $i \in \N$.
\item[4)] The divisor $R^{-}$ is exceptional $/X$ and the codimension of $\Tilde{f}(\text{Supp} R^{-})$ in $\Tilde{Y}$ is $\geq 2$.
\end{itemize}
Now, for a sufficiently divisible $i \in \N$, we have
\[
\Tilde{f}_{*}\mathcal{O}_{\Tilde{X}}(ipK_{\Tilde{X}} +iR^{-}) 
= \Tilde{f}_{*}\mathcal{O}_{\Tilde{X}}(ip\Tilde{f}^{*}(K_{\Tilde{Y}} + B + L) +iR^{+})\]
\begin{equation}\label{eq6}
=\mathcal{O}_{\Tilde{Y}}(ip(K_{\Tilde{Y}} + B + L)).
\end{equation}
On the other hand, by applying \cite[Theorem 5.2]{HSP}, we deduce that 
\begin{equation}\label{eq7}
(h \circ \sigma)_{*}\mathcal{O}_{\Tilde{X}}(ipK_{\Tilde{X}} +iR^{-}) = h_{*}\mathcal{O}_{X}(ipK_{X}) = \mathcal{O}_{A}.
\end{equation}
Combining equations (\ref{eq6}) and (\ref{eq7}), we obtain
\[
(\operatorname{alb_{Y}} \circ \tau)_{*}\mathcal{O}_{\Tilde{Y}}(ip(K_{\Tilde{Y}} + B + L))=  \mathcal{O}_{A}.
\]
If $B+L$ is $\Q$-effective, then 
\[
0= \kappa(X) = \kappa(\Tilde{X}) = \kappa(\Tilde{X}, ipK_{\Tilde{X}} +iR^{-}) = \kappa(\Tilde{Y},ip(K_{\Tilde{Y}} + B + L) )
\]
\[
\geq \kappa(\Tilde{Y}) = \kappa(Y) \geq 0.
\]
For the rest of the proof, we assume that $B+L$ is not $\Q$-effective, in particular, $L$ is nef and not $\Q$-effective. We write 
\[
L = \tau^{*}\tau_{*}L + E_{L}, \hspace{0.2cm}
B= \tau^{*}\tau_{*}B +E_{B},\]
and \[
K_{\Tilde{Y}} = \tau^{*}K_{Y} +E_{Y},
\]
for some exceptional divisors $E_{L}$, $E_{B}$ and $E_{Y}$. For simplicity, we write $\tau_{*}L :=L_{Y}$ and $\tau_{*}B := B_{Y}$. We have $L_{Y}$ is nef and $B_{Y}$ is effective. Thus, from equation (\ref{eq5}), we get
\[
pK_{\Tilde{X}} + R^{-} = p\Tilde{f}^{*}(\tau^{*}(K_{Y} + L_{Y} + B_{Y}) + E_{Y} + E_{L} +E_{B}) + R^{+}.
\]
We set $E: = E_{Y} + E_{L} +E_{B} = E^{+} - E^{-}$, where  $E^{+}$ and $E^{-}$ are the positive and negative parts respectively. Furthermore, we can rewrite the previous equation as follows
\begin{equation}\label{eq8}
pK_{\Tilde{X}} + R^{-} + p\Tilde{f}^{*}E^{-} = p\Tilde{f}^{*}(\tau^{*}(K_{Y} + L_{Y} + B_{Y}) +E^{+}) + R^{+}.
\end{equation}
For a sufficiently divisible $i \in \N$, we have
\begin{equation}\label{eq9}
(h \circ \sigma)_{*}\mathcal{O}_{\Tilde{X}}(ipK_{\Tilde{X}} +iR^{-} + ip\Tilde{f}^{*}E^{-}) = h_{*}\mathcal{O}_{X}(ipK_{X}) = \mathcal{O}_{A}.
\end{equation}
By equations (\ref{eq8}) and (\ref{eq9}), we get
\[
(h \circ \sigma)_{*}\mathcal{O}_{\Tilde{X}}(ip\Tilde{f}^{*}(\tau^{*}(K_{Y} + L_{Y} + B_{Y}) +E^{+}) + iR^{+}) 
\]
\[
= (\operatorname{alb}_{Y} \circ \tau \circ \Tilde{f})_{*}\mathcal{O}_{\Tilde{X}}(ip\Tilde{f}^{*}(\tau^{*}(K_{Y} + L_{Y} + B_{Y}) +E^{+}) + iR^{+}) 
\]
\[
=(\operatorname{alb}_{Y})_{*}\mathcal{O}_{Y}(ip(K_{Y} + L_{Y}+ B_{Y})) = \mathcal{O}_{A}.
\]
It follows that 
\begin{equation}\label{eq10}
\kappa(ip(K_{Y} +L_{Y} + B_{Y})_{|_{G}}) =0,
\end{equation}
where $G$ is a general fiber of $\operatorname{alb}_{Y}$. We know that $K_{Y}$ is semi-ample, and in particular $K_{G}$ is semi-ample as well. Thus, for a sufficiently large $i\in N$, we have $ipK_{G} = 0$. So from equation (\ref{eq10}), we get 
\begin{equation}\label{eq11}
\kappa(ip(L_{Y} + B_{Y})_{|_{G}}) =0.
\end{equation}
We denote ${L_{Y}}_{|_G} = L_{G}$ and ${B_{Y}}_{|_G} = B_{G}$.  Note that $G$ is an elliptic curve or a minimal surface of $\kappa(G)=0$ belongs to one of the four classes: K$3$ surface, abelian surface, Enriques surface, or bi-elliptic surface. 
\begin{itemize}
    \item Case 1: $G$ is a K$3$ surface. Recall that $L_{G}$ is nef and thus $L_{G}$ is semi-ample. Since $B_{G}$ is $\Q$-effective, we have $\kappa(L_{G}) =0$ by equation (\ref{eq11}), and consequently $ipL_{G} = \mathcal{O}_{G}$. Hence $ipL_{Y} = \operatorname{alb}_{Y}^{*}M$ for some non-torsion element $M$ of $\operatorname{Pic}^{0}(A)$.
    \item Case 2: $G$ is an abelian surface. We prove that the nef divisor $L_{G}$ is $\Q$-effective. Indeed, we have $ipB_{G}$ is semi-ample. If $\kappa(B_{G}) =0$, then $ipB_{G} = \mathcal{O}_{G}$, and by the fact that $\kappa(ip(L_{G} + B_{G})) =0$, we deduce that $\kappa(L_{G}) =0$. If $\kappa(B_{G}) = 1$, then the linear system $|ipB_{G}|$ defines a fibered surface $\phi: G \to G^{'}$ with elliptic fibers $G^{''}$, where $G^{'}$ is a an elliptic curve, and $ipB_{G}:= \phi^{*}B_{G^{'}}$ for some ample divisor $B_{G^{'}}$. Hence, combining with equation (\ref{eq11}), we obtain $\kappa({L_{G}}_{|_{G^{''}}}) \geq 0$. If $\kappa({L_{G}}_{|_{G^{''}}}) =0$, then for a sufficiently integer $i\in \N$, we have $ip{L_{G}}_{|_{G^{''}}}= \mathcal{O}_{G^{''}}$. Thus, $ipL_{G} = \phi^{*}L_{G^{'}}$ for a nef divisor $L_{G^{'}}$. Therefore, \[\kappa(ip(L_{G} +B_{G})) = \kappa(\phi^{*}(L_{G^{'}} + B_{G^{'}})) \]\[= \kappa(L_{G^{'}} + B_{G^{'}}) =1,\] which is a contradiction. If $\kappa({L_{G}}_{|_{G^{''}}}) =1$, then $L_{G}$ is nef and $\phi$-big. By applying Fujita type results on positivity of direct image sheaves, we deduce that $\phi_{*}\mathcal{O}_{G}(ipL_{G})$ is a nef rank 2 vector bundle on $G^{'}$, since \[\phi_{*}\mathcal{O}_{G}(ipL_{G}) = \phi_{*}\mathcal{O}_{G}(K_{G/G^{'}}+ (ip-1)K_{G/G^{'}} +ip L_{G}).\] Hence \[
    h^{0}(G, ip(L_{G} +B_{G}))= h^{0}(G^{'}, \phi_{*}\mathcal{O}_{G}(ip(L_{G} +B_{G})))
    \]
    \[
    =h^{0}(G^{'}, \phi_{*}\mathcal{O}_{G}(ipL_{G}) \otimes B_{G^{'}}) \geq \deg (\phi_{*}\mathcal{O}_{G}(ipL_{G}) \otimes B_{G^{'}}) >1,
    \]
    which contradicts equation (\ref{eq11}). If $\kappa(B_{G}) =2$, then $B_{G}$ is big. Hence, $\kappa(B_{G} +L_{G}) =2$ since $L_{G}$ is pseudo-effective, which is a contradiction. We deduce that $L_{G}$ is $\Q$-effective with $\kappa(L_{G})=0$. Thus $ipL_{G}=\mathcal{O}_{G}$. Therefore, $ipL_{Y} ={\operatorname{alb}_{Y}}^{*}M$ for a non-torsion element $M$ of $\operatorname{pic}^{0}(A)$.
\end{itemize}
Finally, in both cases above, we obtain
$ipL_{Y} = \operatorname{alb}_{Y}^{*} M$,
for a non-torsion element \(M \in \operatorname{Pic}^{0}(A)\).
The proof in the case when \(G\) is an Enriques surface is the same as in the case when \(G\) is a K3 surface, and the proof for the case when \(G\) is a bi-elliptic surface or an elliptic curve is the same as when \(G\) is an abelian surface. Hence $ipL_{Y} = \operatorname{alb}_{Y}^{*} M$,
for a non-torsion element \(M \in \operatorname{Pic}^{0}(A)\).
Consequently, we have
\[
ipK_{\Tilde{X}} + iR^{-} + ip\Tilde{f}^{*}E^{-} - \sigma^{*}(h^{*}M)
\]
\[= 
ipK_{\Tilde{X}} + iR^{-} + ip\Tilde{f}^{*}E^{-} - \Tilde{f}^{*}(\tau^{*}({\operatorname{alb}_{Y}}^{*}M)) 
\] 
\[
=ipK_{\Tilde{X}} + iR^{-} + ip\Tilde{f}^{*}E^{-} - ip\Tilde{f}^{*}(\tau^{*}L_{Y})\]
\[= ip\Tilde{f}^{*}(\tau^{*}(K_{Y}  + B_{Y}) +E^{+}) + iR^{+},
\]
where the last equation follows from equation (\ref{eq8}).
Applying \cite[Theorem 3.1]{campanapeternell}, we deduce that
\[
0=\kappa(X) \geq \kappa(ipK_{X}-h^{*}M)
\]
\[
=\kappa(ipK_{\Tilde{X}} + iR^{-} + ip\Tilde{f}^{*}E^{-} - \sigma^{*}(h^{*}M))
\]
\[
=\kappa(ip\Tilde{f}^{*}(\tau^{*}(K_{Y}  + B_{Y}) +E^{+}) + iR^{+})
\]
\[
=\kappa(ip(\tau^{*}(K_{Y}  + B_{Y}) +E^{+}))
\]
\[
\kappa(K_{Y} + B_{Y}) \geq \kappa(Y) \geq 0.
\]
Finally, we obtain $\kappa(Y) =0$.
\par Step 5. Induction: We proceed by induction on the dimension of \( X \). 
If \( \dim X = 1 \), then there is nothing to prove. 
Now assume that Theorem~\ref{maintheorem} holds in all dimensions \( < n \). Let \( f: X \to Y \) be the given algebraic fiber space with \( \dim X = n \). 
If \( \kappa(X) = -\infty \), then by inequality~(\ref{eq3}) we have \( \kappa(F) = -\infty \), and the theorem follows from Step~1. 
If \( \kappa(X) = 0 \), then by Step~4 we deduce that \( \kappa(Y) = 0 \), and by Step~3 the theorem is proved. 
Now assume \( \kappa(X) > 0 \). Let $g: X \to Z$ be the Iitaka fibration of $X$, and let $X_{Z}$ be the general fiber of $g$. We know that $\kappa(X_{Z}) = 0$ and $\kappa(X) = \dim Z > 0$.
\[
\begin{tikzcd}
X \arrow[r, "g"] \arrow[d, "f"] & Z  \\
 Y & 
\end{tikzcd}
\]
We define the restriction $f_{|_{X_{Z}}} : X_{Z} \to B :=f(X_{Z})$. 
\begin{itemize}
\item $B$ is a point. In this case, $X_{Z}$ is contracted by $f$. This means that we can define a map $\psi: Z \to Y$, which can be assumed to be regular after resolving the indeterminacy locus. Note that  $g(F)$ is a general fiber of $\psi$. Hence, we define the restriction $g_{|_{F}}: F \to g(F)$ with general fiber $X_{Z}$. It follows that 
\begin{equation}\label{eq12}
\kappa(F) \leq \kappa(X_{Z}) + \dim g(F) = \dim g(F).
\end{equation}
However, 
\begin{equation}\label{eq13}
\kappa(X) =\dim Z = \dim g(F) + \dim Y.
\end{equation}
Combining equations (\ref{eq12}) and (\ref{eq13}), we deduce the desired inequality
\[
\kappa(X) \geq \kappa(F) + \dim Y \geq \kappa(F) + \kappa(Y).
\]
    \item $B$ is not a point. By Lemma \ref{lemma}, there exist a general subvariety $\Gamma$ of $Z$ , a projective variety $\widehat{Y}$ and morphisms $\gamma: \widehat{Y} \to \Gamma  , \beta : \widehat{Y} \to Y$ such that
$\gamma$ is a surjective morphism with connected fibers, and $\gamma^{-1}(z) \simeq f(g^{-1}(z)) =B$ for $z \in \Gamma$, and $\beta$ is generically finite. By Proposition \ref{propo}, $\alpha_{f}(B) \leq \alpha_{f}(\widehat{Y}) \leq \alpha_{f}(Y) \leq 2$. Thus applying the induction hypothesis to $f_{|_{X_{Z}}}$, and obtain
\begin{equation}\label{eq14}
0 = \kappa(X_{Z}) \geq \kappa(X_{Y,Z}) + \kappa(B).
\end{equation}
Here $X_{Y,Z}$ is a general fiber of $f_{|_{X_{Z}}}$. Clearly $\kappa(X_{Y,Z}) \geq 0$. Next, we will see that $\kappa(B) \geq 0$. Indeed, by Mori's easy addition formula, we get
\begin{equation}\label{eq15}
0 \leq \kappa(Y) \leq \kappa(\widehat{Y}) \leq \kappa(B) + \dim \Gamma.
\end{equation}
This implies $\kappa(B) \geq 0$. By inequality (\ref{eq14}), we deduce that 
\begin{equation}\label{eq16}
\kappa(X_{Y,Z}) = \kappa(B) = 0.
\end{equation}
Now, we consider the restriction $g_{|_F}: F \to g(F)$ with general fiber $X_{Y,Z}$. Therefore, we apply Mori's easy addition formula to $g_{|_{F}}$, we get
\[
\kappa(F) \leq \kappa(X_{Y,Z}) +\dim g(F).
\]
Then applying equation (\ref{eq16}), we obtain
\begin{equation}\label{eq17}
\kappa(F) \leq \dim g(F).\end{equation}

Furthermore, we have

\[
\dim g(F) = \dim F - \dim X_{Y,Z}
\]
\[
= \dim X - \dim Y - (\dim X_{Z} -\dim B)
\]
\[
=\dim Z + \dim B - \dim Y
\]
\[
=\kappa(X) + \dim B - \dim Y.
\]
Combining these previous equations with (\ref{eq17}), we get
\[
\kappa(X) \geq \kappa(F) +\dim Y - \dim B \]
\[= \kappa(F) + \dim \Gamma \geq \kappa(F) + \kappa(Y).
\]
The last inequality follows from inequalities (\ref{eq15}) and (\ref{eq16}).
\end{itemize}
\end{proof}
\begin{remark}
We point out that Theorem $\ref{schnellthm}$ is useful for the Iitaka Conjecture \ref{Conj1}. Particularly in our situation. For the reader's interest, we explain this fact in the following corollary.
\end{remark}
\begin{corollary}
Theorem \ref{schnellthm} implies that the Iitaka conjecture $C_{n,m}$ holds when $\alpha(Y)\geq m-2$, provided that the Iitaka variety of $Y$ is not uniruled.
\end{corollary}
\begin{proof}
The proof is by induction. Steps $1, 2$, and $3$ are the same as in Theorem \ref{maintheorem}. 
\begin{itemize}
    \item Step 4. $\kappa(Y) >0 \implies \kappa(X) \geq \kappa(F) + \kappa(Y)$. Indeed, by \cite{birkarchen} and \cite{Huzhengyu}, we can run the minimal model program (MMP) for $Y$ and we obtain a good MMP $Y^{'}$ of $Y$.  Without loss of generality, we write $Y^{'}= Y$, so that $K_{Y}$ is semi-ample. We then take the Iitaka fibration $I$ of $Y$, and we have the following commutative diagram.
    \[
\begin{tikzcd}
X \arrow[r, "f"] \arrow[rd, "J"] & Y \arrow[d, "I"] \\
 & I(Y)
\end{tikzcd}
\]
We may assume that $I$ is an algebraic fiber space, where $I(Y)$ is the Iitaka variety of $Y$. Thus, for some integer $m \in \N$, we have,
$mK_{Y} = I^{*}M$, for an ample divisor $M$ on $I(Y)$. Consequently, \[
mK_{X/Y} = mK_{X} - mf^{*}K_{Y} = mK_{X} - J^{*}M.
\]
By Step 1, we have $\kappa(F) \geq 0$. Thus $K_{X/Y}$ is pseudo-effective. 
Since $I(Y)$ is not uniruled by hypothesis,  we apply Theorem \ref{maintheorem} to deduce that $mK_{X} - J^{*}M$ is effective for a sufficiently large $m \in \N$. Applying \cite[Proposition 1.14]{Mori}, we deduce that 
\begin{equation}\label{eq18}
\kappa(X) = \kappa(X_{I(Y)}) + \dim I(Y) = \kappa(X_{I(Y)}) + \kappa(Y)
\end{equation}
where $X_{I(Y)}$ is a general fiber of $J$. We define the restriction $f_{|_{X_{I(Y)}}}: X_{I(Y)} \to Y_{I(Y)}$ with general fiber $F$. Here, $Y_{I(Y)}$ is a general fiber of $I$ with $\kappa(Y_{I(Y)}) =0$, and it satisfies the same Albanese dimension property as $Y$ by Proposition \ref{propo}. Therefore, by induction, we have
\begin{equation}\label{eq19}
\kappa(X_{I(Y)}) \geq \kappa(F) + \kappa(Y_{I(Y)}) = \kappa(F).
\end{equation}
Finally, by combining (\ref{eq18}) and (\ref{eq19}), we obtain the desired result
\[
k(X) \geq k(F) + k(Y).
\]
\end{itemize}
\end{proof}
\bibliographystyle{plain}
\bibliography{slope paper}
\end{document}